\documentclass{article}
 \usepackage{amsmath}

\newtheorem{thm}{Theorem}[section]

\newtheorem{lem}[thm]{Lemma}

\newcommand{\be}{\begin{equation}}
\newcommand{\ee}{\end{equation}}
\newcommand{\ben}{\begin{enumerate}}
\newcommand{\een}{\end{enumerate}}
\newcommand{\beq}{\begin{eqnarray}}
\newcommand{\eeq}{\end{eqnarray}}
\newcommand{\beqn}{\begin{eqnarray*}}
\newcommand{\eeqn}{\end{eqnarray*}}

\newcommand{\pa}{\partial}

\newcommand{\pxi}{ {\pa \over \pa x^i}}

\newcommand{\qed}{\hspace*{\fill}Q.E.D.}  

\begin{document}
\title{On a Class of Two-Dimensional Einstein Finsler
Metrics of Vanishing S-Curvature}
\author{Guojun Yang  }
\date{}
\maketitle
\begin{abstract}
   An $(\alpha,\beta)$-metric is defined by a Riemannian metric
$\alpha$ and $1$-form $\beta$.  In this paper, we study a known
class of two-dimensional $(\alpha,\beta)$-metrics of vanishing
S-curvature. We determine the local structure of those metrics and
show that those metrics are Einsteinian (equivalently, isotropic
flag curvature) but generally are not Ricci-flat.

\

{\bf Keywords:}  $(\alpha,\beta)$-Metric, Einstein Metric,
 S-Curvature, Flag Curvature

 {\bf MR(2000) subject classification: }
53B40
\end{abstract}

\section{Introduction}

In Finsler geometry, Einstein metrics are defined in a natural way
as that in Riemann geometry. An $n$-dimensional Finsler metric $F$
is called an Einstein metric if its Ricc curvature $Ric$ is
isotropic,
 $$
 Ric=(n-1)\lambda F^2,
 $$
where $\lambda=\lambda(x)$ is a scalar function. $F$ is called of
Ricci constant, if $\lambda=constant$. In particular, $F$ is
called Ricci-flat if $\lambda=0$. It is well known that in
dimension $n\ge 3$, every Einstein Riemann metric is of Ricci
constant, and  every 3-dimensional Einstein Riemann metric is of
constant sectional curvature. We do not know whether it is still
true for any Finsler metrics. It has been shown that many Finsler
metrics have such a similar property as Riemann metrics, among
which, two important cases are Randers metrics (cf. \cite{BR2})
and square metrics (cf. \cite{CSZ}). An Einstein square metric in
$n\ge 2$ is always Ricci-flat (cf. \cite{CSZ} \cite{CST}), but it
is not necessarily the case for a Randers metric (cf. \cite{BR2}
\cite{BRS}).

The S-curvature is originally introduced for the volume comparison
theorem (\cite{shen1}), and it is a non-Riemannian quantity which
plays an important role in Finsler geometry (cf. \cite{CS1}
\cite{CMS} \cite{shen1}--\cite{Y2}). For a Finsler manifold,  the
flag curvature
 is an analogue of sectional curvature for a Riemann manifold.
The flag curvature and the S-curvature are closely related. It is
proved that, for a Finsler
 manifold $(M,F)$ of scalar flag curvature, if $F$ is of  isotropic
 S-curvature ${\bf S}=(n+1)c(x)F$ for a scalar function $c(x)$
 on $M$, then the flag curvature must be in the following form
 \be\label{y1}
 {\bf K}=\frac{3c_{x^m}y^m}{F}+\tau(x),
 \ee
where $\tau(x)$ is a scalar function on $M$ (\cite{CMS}). Clearly,
if $F$ is of constant S-curvature and of scalar flag curvature,
then by (\ref{y1}), ${\bf K}=\tau(x)$ is isotropic and ${\bf
K}=constant$ in $n\ge 3$.

An  $(\alpha,\beta)$-metric is defined by a Riemannian metric
 $\alpha=\sqrt{a_{ij}(x)y^iy^j}$ and a $1$-form $\beta=b_i(x)y^i\ne 0$ on a manifold
 $M$, which can be expressed in the following form:
 $$F=\alpha \phi(s),\ \ s=\beta/\alpha,$$
where $\phi(s)$ is a function satisfying certain conditions such
that $F$ is positive definite on $TM-0$ (see \cite{shen2}). Some
recent studies show that many Einstein $(\alpha,\beta)$-metrics
are Ricci-flat.  In \cite{CST}, it proves that an Einstein
$(\alpha,\beta)$-metric with $\phi(s)$ being a non-linear
polynomial must be Ricci-flat. By this result, it is ever believed
that any Einstein $(\alpha,\beta)$-metrics with $\phi(s)$ being
non-linear analytic must be Ricci-flat. But this is not true. In
\cite{Y3}, the present author studies a class of
$(\alpha,\beta)$-metric $F$ with $\phi(s)=(1+s)^p$, where $p\ne 0$
is a real number, and shows that  a two-dimensional Einstein
square-root metric $F=\sqrt{\alpha(\alpha+\beta)}$ ($p=1/2$) is
generally not Ricci-flat (also see the following Theorem
\ref{th1}).

 Generally, two-dimensional Finsler metrics have some different special
curvature properties from higher dimensions (cf. \cite{Y2}
\cite{Y4}--\cite{Y6}). By definition, every 2-dimensional Finsler
metric is of scalar flag curvature ${\bf K}={\bf K}(x,y)$, but
generally {\bf K} is not isotropic. A two-dimensional metric $F$
is an Einstein metric if and only if $F$ is of isotropic flag
curvature. By (\ref{y1}), if $F$ is of constant S-curvature, then
{\bf K} is isotropic. Conversely, if a Randers metric
$F=\alpha+\beta$ is of isotropic flag curvature, then $F$ is of
constant S-curvature (\cite{BR2}). In \cite{Y}, we
 construct a family of two-dimensional Randers metrics which
are of isotropic flag curvature ${\bf K}={\bf K}(x)\ne constant$.
In \cite{Y3}, we prove that a 2-dimensional square-root metric
$F=\sqrt{\alpha(\alpha+\beta)}$ is of isotropic flag curvature if
and only if $F$ is of vanishing S-curvature. In \cite{Y2}, we
investigate again the known characterization (\cite{CS1}) for
$(\alpha,\beta)$-metrics of isotropic S-curvature, and obtain one
more class of two-dimensional $(\alpha,\beta)$-metrics of
vanishing S-curvature with $\phi(s)$ defined by (\ref{y2}) below.
In this paper, we will study certain curvature properties  of
such a class of $(\alpha,\beta)$-metrics in the following theorem.

\begin{thm}\label{th1}
 Let $F=\alpha \phi(s)$, $s=\beta/\alpha$, be a two-dimensional
  $(\alpha,\beta)$-metric with $\phi(s)$ being defined by
  \be\label{y2}
      \phi(s)=\big\{(1+k_1s^2)(1+k_2s^2)\big\}^{\frac{1}{4}}e^{\int^s_0\tau(s)ds},
     \ee
     where
      $$
       \tau(s):=\frac{\pm\sqrt{k_2-k_1}}
      {2(1+k_1s^2)\sqrt{1+k_2s^2}},
      $$
     and $k_1$ and $k_2$ are constants with $k_2>k_1$. Suppose
     $F$ is of isotropic S-curvature. Then locally we have
     \beq\label{y3}
  \alpha&=&\frac{\sqrt{B}}{(1+k_1B)^{\frac{3}{4}}
  (1+k_2B)^{\frac{1}{4}}}\sqrt{\frac{(y^1)^2+(y^2)^2}{u^2+v^2}},\label{y3}\\
  \beta&=&\frac{B}{(1+k_1B)^{\frac{3}{4}}
  (1+k_2B)^{\frac{1}{4}}}\frac{uy^1+vy^2}{u^2+v^2},\label{y4}
  \eeq
  where
    $B=B(x),u=u(x),v=v(x)$ are some scalar functions
    satisfying
   \be\label{y5}
  u_1=v_2,\ \ u_2=-v_1,\ \ uB_1+vB_2=0,
   \ee
where $u_i:=u_{x^i},v_i:=v_{x^i}$ and $B_i:=B_{x^i}$. Further, the
S-curvature ${\bf S}=0$ and the  flag curvature ${\bf K}$ is
isotropic given by
 \be\label{y6}
 {\bf
 K}=-\frac{(u^2+v^2)\sqrt{1+k_2B}}{4B^2}\Big\{2\sqrt{1+k_1B}(B_{11}+B_{22})
 -\frac{(u^2+v^2)(2+3k_1B)}{B\sqrt{1+k_1B}}\big(\frac{B_1}{v}\big)^2\Big\},
 \ee
 where $B_{ij}:=B_{x^ix^j}$.
\end{thm}

Theorem \ref{th1} gives a class of two-dimensional Einstein
Finsler metrics, but generally they are not Ricci-flat. We can
easily find $u,v,B$ satisfying (\ref{y5}), for example, $u=-x^2,\
v=x^1$ and $B=(x^1)^2+(x^2)^2$.

If $v=0$ in (\ref{y6}), then $B_1/v$ can be replaced by $-B_2/u$
since (\ref{y5}).  Further, $F$ is positively definite if and only
if $1+k_1b^2>0$ (Lemma \ref{lem31} below). If
$B(=||\beta||^2_{\alpha})=constant$, then $\alpha$ is flat and
$\beta$ is parallel with respect to $\alpha$ (Lemma \ref{lem32}
below).

In \cite{Y2}, we prove that $F$ in Theorem \ref{th1} is of
isotropic S-curvature if and only if $F$ is of vanishing
S-curvature. Take $k_1=-1,k_2=0$, then $F$ in Theorem \ref{th1}
becomes $F=\sqrt{\alpha(\alpha+\beta)}$, which is called a
square-root metric (\cite{Y3}). In \cite{Y3}, when we study a
class of $(\alpha,\beta)$-metrics of Einstein-reversibility,
Theorem \ref{th1} has been actually proved for
$F=\sqrt{\alpha(\alpha+\beta)}$, and we also show that the
converse is also true for square-root metrics, namely,  Einstein
square-root metrics must be of vanishing S-curvature.

The converse of Theorem \ref{th1} might be true. We have a
conjecture: {\it if the two-dimensional $(\alpha,\beta)$-metric
$F$ defined by (\ref{y2}) is Einsteinian, then $F$ must be of
vanishing S-curvature.} Some special cases can be verified, for
example, $k_1=-1,k_2=0$ as shown above. We have also verified
another more complicated  special case: $k_1=0,k_2=4$. In general
case, we only need to prove (\ref{r00}) below holds if $F$ defined
by (\ref{y2}) is an Einstein metric.

\section{Preliminaries}

For a Finsler metric $F$, the Riemann curvature $R_y=R^i_{\
k}(y)\frac{\pa}{\pa x^i}\otimes dx^k$ is defined by
 $$
 R^i_{\ k}:=2\frac{\pa G^i}{\pa x^k}-y^j\frac{\pa^2G^i}{\pa x^j\pa
 y^k}+2G^j\frac{\pa^2G^i}{\pa y^j\pa y^k}-\frac{\pa G^i}{\pa y^j}\frac{\pa G^j}{\pa
 y^k},
 $$
 where $G^i$ are called the geodesic coefficients as follows
 $$
  G^i:=\frac{1}{4}g^{il}\big\{[F^2]_{x^ky^l}y^k-[F^2]_{x^l}\big\}.
  $$
Then the Ricci curvature ${\bf Ric}$ is defined by ${\bf
Ric}:=R^k_{\ k}$. A Finsler metric is called of scalar flag
curvature if there is a function ${\bf K}={\bf K}(x,y)$ such that
 \be\label{yR}
  R^i_{\ k}={\bf K}F^2(\delta^i_k-F^{-2}y^iy_k), \ \ y_k:=(1/2F^2)_{y^iy^k}y^i.
  \ee
In two-dimensional case, we have
 \be\label{K}
{\bf K}=\frac{{\bf Ric}}{F^2}.
 \ee

Under the Hausdorff-Busemann volume form $dV=\sigma_F(x)dx^1\wedge
... \wedge dx^n$, where
  $$\sigma_F(x):=\frac{Vol(B^n)}{Vol\big\{(y^i)\in
  R^n|F(y^i\pxi|_x)<1\big \}},$$
 the S-curvature is defined by
$${\bf S}:=\frac{\pa G^m}{\pa y^m}-y^m\frac{\pa}{\pa
x^m}(ln\sigma_F).$$
 ${\bf S}$ is said to be isotropic if there is a scalar function
 $c(x)$ on $M$ such that
 $${\bf S}=(n+1)c(x)F.$$
 If $c(x)$ is a constant, then  $F$ is called of  constant
 S-curvature.

For a Riemannian  $\alpha =\sqrt{a_{ij}y^iy^j}$ and a $1$-form
$\beta = b_i y^i $, let
 $$r_{ij}:=\frac{1}{2}(b_{i|j}+b_{j|i}),\ \ s_{ij}:=\frac{1}{2}(b_{i|j}-b_{j|i}),\ \
 r^i_{\ j}:=a^{ik}r_{kj},\ \  s^i_{\ j}:=a^{ik}s_{kj},$$
 $$q_{ij}:=r_{im}s^m_{\ j}, \ \ t_{ij}:=s_{im}s^m_{\ j},\ \
 r_j:=b^ir_{ij},\ \  s_j:=b^is_{ij},$$
 $$
  q_j:=b^iq_{ij}, \ \ r_j:=b^ir_{ij},\ \ t_j:=b^it_{ij},
 $$
 where we define $b^i:=a^{ij}b_j$, $(a^{ij})$ is the inverse of
 $(a_{ij})$, and $\nabla \beta = b_{i|j} y^i dx^j$  denotes the covariant
derivatives of $\beta$ with respect to $\alpha$. Here are some of
our conventions in the whole paper. For a general tensor $T_{ij}$
as an example, we  define $T_{i0}:=T_{ij}y^j$ and
$T_{00}:=T_{ij}y^iy^j$, etc. We use $a_{ij}$ to raise or lower the
indices of a tensor.

\begin{lem}(\cite{BR2})\label{lem21}
By Ricci identities we have
 $$s_{ij|k}=r_{ik|j}-r_{jk|i}-b^l\bar{R}_{klij},$$
 $$s^k_{\ 0|k}=r^k_{\ k|0}-r^k_{\ 0|k}+b^l\bar{R}ic_{l0},$$
 $$b^ks_{0|k}=r_ks^k_{\ 0}-t_0+b^kb^lr_{kl|0}-b^kb^lr_{k0|l},$$
 $$s^k_{\ |k}=r^k_{\ |k}-t^k_{\ k}-r^i_{\ j}r^j_{\ i}-b^ir^k_{\
 k|i}-b^kb^i\bar{R}ic_{ik}.$$
 where $\bar{R}$ denotes the Riemann curvature tensor of $\alpha$.
\end{lem}

\begin{lem}(\cite{Y2})\label{lem22}
Let $F=\alpha \phi(s)$, $s=\beta/\alpha$, be a two-dimensional
  $(\alpha,\beta)$-metric on a manifold $M$ with $\phi(s)$ being defined by
  (\ref{y2}). Then $F$ is of isotropic S-curvature ${\bf S}=3c(x)F$ if and only if
 \be\label{r00}
 r_{ij}=\frac{3k_1+k_2+4k_1k_2b^2}{4+(k_1+3k_2)b^2}(b_is_j+b_js_i).
 \ee
 In this case, ${\bf S}=0$
\end{lem}

\begin{lem}(\cite{Y})\label{lem23}
 Let $\alpha=\sqrt{a_{ij}y^iy^j}$ be an $n(\ge 2)$-dimensional Riemann metric which is locally conformally
 flat. Locally we express $a_{ij}=e^{2\sigma(x)}\delta_{ij}$. Then $W_0=W_iy^i$ is a conformal 1-form of $\alpha$
 satisfying
  $$W_{0|0}=-2c\alpha^2,$$
where $c=c(x)$ is a scalar function and the covariant derivative
is taken with respect to the Levi-Civita connection of $\alpha$,
if and only if
  \be\label{W}
  \frac{\pa W^i}{\pa x^j}+\frac{\pa W^j}{\pa
  x^i}=0 \quad{\rm (} \forall \thinspace i\ne  j{\rm )},\ \ \  \frac{\pa W^i}{\pa x^i}=\frac{\pa W^j}{\pa
  x^j} \quad{\rm (} \forall \thinspace i, j{\rm )},
  \ee
  where $W^i:=a^{ij}W_j$. In this case, $c$ is given by
   \be\label{c}
    c(x)=-\frac{1}{2}\big[\tau(x)+W^r\sigma_r\big], \ \ \ (\tau:=\frac{\pa W^1}{\pa x^1},\ \ \sigma_i:=\sigma_{x^i}).
    \ee
\end{lem}

\section{Proof of Theorem \ref{th1}}

We first give a lemma to show the positively definite condition of
the $(\alpha,\beta)$-metric defined by (\ref{y2}).

\begin{lem}\label{lem31}
Let $F=\alpha \phi(s)$, $s=\beta/\alpha$, be a two-dimensional
  $(\alpha,\beta)$-metric on a manifold $M$ with $\phi(s)$ being defined by
  (\ref{y2}). Then $F$ is positively definite  on $TM-0$ if and only if
 \be\label{pdc}
 1+k_1b^2>0.
 \ee
\end{lem}

{\it Proof :} In \cite{shen2}, it is shown that an
$(\alpha,\beta)$-metric $F=\alpha \phi(\beta/\alpha)$ is
positively definite if and only if for $|s|\le b$,
  \be\label{y11}
 \phi(s)>0,\ \ \ \phi(s)-s\phi'(s)>0,\ \ \
 \phi(s)-s\phi'(s)+(b^2-s^2)\phi''(s)>0.
  \ee
Now let $\phi(s)$ be given by (\ref{y2}) in (\ref{y11}). Clearly
we have $\phi(s)>0$ for $|s|\le b$. It is easy to see that
 $$F \ \text{is defined on} \ TM-0 \Longleftrightarrow 1+k_1s^2>0
\Longleftrightarrow (\ref{pdc}).
 $$
In the following, we assume (\ref{pdc}). Then we only need to
prove the second and the third inequalities hold.

 The second
inequality in (\ref{y11}) is equivalent to
 $$\big(\sqrt{1+k_2s^2}\mp \sqrt{k_2-k_1}s\big)\big(\sqrt{1+k_2s^2}\pm
 \frac{\sqrt{k_2-k_1}}{2}s\big)>0.$$
Clearly, the above inequality holds by (\ref{pdc}).

The third inequality in (\ref{y11}) is equivalent to
 \be\label{y12}
 f(s^2)\sqrt{1+k_2s^2}+g(s)>0, \ \ (|s|\le b),
 \ee
where
 $$f(s^2):=\big[2k_1k_2(k_1+k_2)b^2+3k_1^2+2k_2^2-k_1k_2\big]s^4+\big[k_1(9k_2-k_1)b^2+3k_2+5k_1\big]s^2+(3k_2+k_1)b^2+4,$$
 $$g(s):=-4a_1(1+k_1b^2)(1+k_2s^2)^2s,\ \ \
 a_1:=\pm\frac{\sqrt{k_2-k_1}}{2}.$$
To prove (\ref{y12}), we first show $f(s^2)>0$ for $|s|\le b$. Put
 $$f(t)=\widetilde{A}t^2+\widetilde{B}t+\widetilde{C},\ \ \ (0\le t\le b^2),$$
where
 $$
\widetilde{A}:=2k_1k_2(k_1+k_2)b^2+3k_1^2+2k_2^2-k_1k_2,\
\widetilde{B}:=k_1(9k_2-k_1)b^2+3k_2+5k_1,\
\widetilde{C}:=(3k_2+k_1)b^2+4.
 $$
Then by (\ref{pdc}) we have
 $$f(0)=(3k_2+k_1)b^2+4>0,\ \
 f(b^2)=2(1+k_1b^2)(1+k_2b^2)(2+k_1b^2+k_2b^2)>0.
 $$
If $\widetilde{A}= 0$, then the above shows $f(s^2)>0$ for $|s|\le
b$. If $\widetilde{A}\ne 0$, then we get
 $$
 f(-\frac{\widetilde{B}}{2\widetilde{A}})=\frac{a_1^4(1+k_1b^2)\big[(24k_2-k_1)b^2+23\big]}{\widetilde{A}}.
 $$
The above shows $\widetilde{A}> 0$ since $f(0)>0$ or $f(b^2)>0$.
 Thus $f(-\frac{\widetilde{B}}{2\widetilde{A}})>0$. Therefore, $f(s^2)>0$ for
$|s|\le b$ in this case.

Now we prove (\ref{y12}). Without loss of generality, we assume
$a_1>0$ and $0\le s\le b$. Then (\ref{y12}) is equivalent to
$(1+k_2s^2)f^2(s^2)-g^2(s)>0$, which is equivalent to
 $$
 h(t):=\widehat{A}t^2+\widehat{B}t+\widehat{C}>0, \ \ \ (0\le t:=s^2\le b^2),
 $$
where
 $$
\widehat{A}:=4k_1k_2^2(3k_2+k_1)b^4+4k_2(3k_2^2+3k_1^2+2k_1k_2)b^2-6k_1k_2+9k_1^2+13k_2^2,
 $$
 $$
 \widehat{B}:=4k_1k_2(9k_2-k_1)b^4+2(9k_2-k_1)(3k_1+k_2)b^2+12k_1+20k_2,\
 \ \ \widehat{C}:=(4+k_1b^2+3k_2b^2)^2.
 $$
Likewise, we easily get $h(0)>0,\ h(b^2)>0$. Then in a similar way
to the proof of $f(t)>0$ for $0\le t\le b^2$ above, we obtain
$h(t)>0$ for $0\le t\le b^2$.  \qed

\begin{lem}\label{lem32}
 Under the condition of Lemma \ref{lem22}, if there is a
 neighborhood $U$ such that $b=constant$, then $\alpha$ is flat
 and $\beta$ is parallel with respect to $\alpha$ in $U$. In particular,
 if $2(1-k_1k_2b^4)+(k_2-k_1)b^2=
 0$ in $U$, then $\alpha$ is flat
 and $\beta$ is parallel with respect to $\alpha$ in $U$.
\end{lem}

{\it Proof :} If there is a
 neighborhood $U$ such that $2(1-k_1k_2b^4)+(k_2-k_1)b^2=
 0$ in $U$, then $b=constant$ in $U$. Since $b=constant$ is equivalent to $r_i+s_i=0$, by
 (\ref{r00}), we get $b=constant$ is equivalent to
  $$(1+k_1b^2)(1+k_2b^2)s_i=0.$$
  So by the positive definiteness of $F$ (Lemma \ref{lem31}), we
  have $s_i=0$. Since $n=2$, we have (\ref{y26}) below. So
  $s_{ij}=0$, namely, $\beta$ is closed. Then by (\ref{r00}),
  $\beta$ is parallel with respect to $\alpha$, and thus $\alpha$
  is flat.
     \qed

\

By Lemma \ref{lem32} and continuity, we only need to have a
 discussion at those points $x\in M$ with
$2(1-k_1k_2b^4)+(k_2-k_1)b^2\ne 0$.

Now we determine the local structure in Theorem \ref{th1} given by
(\ref{y3})--(\ref{y5}). Since $F$ in Theorem \ref{th1} is of
isotropic S-curvature, $\beta$ satisfies (\ref{r00}) by Lemma
\ref{lem22}. Define a Riemannian metric $\widetilde{\alpha}$ and
1-form $\widetilde{\beta}$ by
 \be\label{y15}
\widetilde{\alpha}:=\alpha, \ \
\widetilde{\beta}:=(1+k_1b^2)^{-\frac{3}{4}}(1+k_2b^2)^{-\frac{1}{4}}\beta.
 \ee
Under the deformation (\ref{y15}), a direct computation shows that
(\ref{r00}) is reduced to
 \be\label{y16}
 \widetilde{r}_{ij}=0.
 \ee
So $\widetilde{\beta}$ is a Killing form with respect to $\alpha$.
Locally  we can express $\alpha$  as
 \be\label{y17}
\alpha:=e^{\sigma}\sqrt{(y^1)^2+(y^2)^2},
 \ee
where $\sigma=\sigma(x)$ is a scalar function. Then by (\ref{W})
in Lemma \ref{lem23}, we have
 \be\label{y18}
 \widetilde{\beta}=\widetilde{b}_1y^1+\widetilde{b}_2y^2=e^{2\sigma}(uy^1+vy^2),
 \ee
where $u=u(x),v=v(x)$ are a pair of scalar functions such that
 $$f(z)=u+iv, \ \ z=x^1+ix^2$$
 is a complex analytic function, and further it follows from (\ref{c}) in Lemma \ref{lem23} and
  (\ref{y16}) that
 $u,\ v$ and $\sigma$ satisfy the following PDEs:
  \be\label{y19}
 u_1=v_2,\ \ \ \ u_2=-v_1,\ \ \ \ u_1+u\sigma_1+v\sigma_2=0.
  \ee
Actually $\sigma$ can be determined  in terms of the triple
$(B,u,v)$, where $B:=b^2$.
 Firstly by (\ref{y15}) and then by (\ref{y17}) and
(\ref{y18}) we get
 \be\label{y20}
||\widetilde{\beta}||^2_{\alpha}=\frac{B}{(1+k_1B)^{\frac{3}{2}}\sqrt{1+k_2B}},\
\ \ \ ||\widetilde{\beta}||^2_{\alpha}=e^{2\sigma}(u^2+v^2).
 \ee
Therefore, by (\ref{y20}) we get
 \be\label{y21}
 e^{2\sigma}=\frac{B}{(u^2+v^2)(1+k_1B)^{\frac{3}{2}}\sqrt{1+k_2B}}.
 \ee
Now it follows from (\ref{y15}), (\ref{y17}), (\ref{y18}) and
(\ref{y21}) that (\ref{y3}) and (\ref{y4})  hold. By (\ref{y21}),
we have
 \be\label{y22}
 \sigma=ln\frac{\sqrt{B}}{\sqrt{u^2+v^2}(1+k_1B)^{\frac{3}{4}}(1+k_2B)^{\frac{1}{4}}}.
 \ee
Plugging (\ref{y22}) into $u_1+u\sigma_1+v\sigma_2=0$ and using
 $u_1=v_2,\ u_2=-v_1$ in (\ref{y19}),
we can write $u_1+u\sigma_1+v\sigma_2=0$ equivalently as
 \be\label{y23}
 \big[2(1-k_1k_2b^4)+(k_2-k_1)b^2\big](uB_1+vB_2)=0.
 \ee
By Lemma \ref{lem32} and continuity, (\ref{y23}) implies
 $$uB_1+vB_2=0.$$
Thus (\ref{y5}) holds.

Next we will use (\ref{K}) to prove that the flag curvature of $F$
in Theorem \ref{th1} is given by (\ref{y6}). The expression of the
Ricci curvature {\bf Ric} for an $(\alpha,\beta)$-metric is
generally very long (see \cite{CST}). We will not write out the
expression of the Ricci curvature of $F$ in Theorem \ref{th1}. To
compute {\bf Ric} of $F$, we need to show the following
quantities:
 \beq\label{y024}
 &&r_{00},\ r_{00|0},\ r_0,\ r_{0|0},\ r^m_{\ m},\ r,\ b^kr_{00|k},\  b^ks_{0|k},\
 b^kq_{0k},\ q_{00},\ q_0,\nonumber\\
  &&\ \ \ \ \ \ \ \ \ \ t_{00},\ t_0,\ t^m_{\ m},\ s_0^2,\
 s_ms^m,\ s^m_{\ 0|m},\ s_{0|0}, \ Ric_{\alpha}.
 \eeq

\begin{lem}\label{lem33}
 In dimension $n=2$, for a pair $\alpha=\sqrt{a_{ij}y^iy^j}$ and
 $\beta=b_iy^i$, there is a $\theta=\theta(x)$ such that
 \be\label{y24}
  s_0^2=\theta(b^2\alpha^2-\beta^2),\ \ \ t_{00}=-\theta \alpha^2.
  \ee
  Further, by  (\ref{y24}), we easily get
   \be\label{y25}
   s_ms^m=\theta b^2,\ \ t_0=-\theta\beta,\ \ t^m_{\ m}=-2\theta.
   \ee
\end{lem}

{\it Proof :} Fix a point $x\in M$ and take  an orthonormal basis
  $\{e_i\}$ at $x$ such that
   $$\alpha=\sqrt{(y^1)^2+(y^2)^2},\ \ \beta=by^1.$$
Then
$$s_0=s_2y^2,\ \ \ b^2\alpha^2-\beta^2=b^2(y^2)^2.$$
So for some $\theta=\theta(x)$, the first formula in (\ref{y24})
holds. Similarly, since $n=2$, it can be easily verified that
 \be\label{y26}
 s_{ij}=\frac{b_is_j-b_js_i}{b^2}.
 \ee
Therefore, by (\ref{y26}),  we have
 \be\label{y27}
 t_{00}=-\frac{s_ms^m\beta^2+b^2s_0^2}{b^4}.
 \ee
Then plugging $s_0^2$ in (\ref{y24}) and $s_ms^m$ in (\ref{y25})
into (\ref{y27}) yields $t_{00}$ in (\ref{y24}).   \qed

\

By (\ref{r00}), we can compute the following quantities:
 $$
r_{00},\ r_{00|0},\ r_0,\ r_{0|0},\ r^m_{\ m},\ r,\ b^kr_{00|k},\
 b^kq_{0k},\ q_{00},\ q_0
 $$
To name a few, we have
 $$
r^m_{\ m}=0,\ \ r=0,\ \
 r_0=\frac{(3k_1+k_2+4k_1k_2b^2)b^2}{4+(k_1+3k_2)b^2}s_0,
 $$
 \beqn
 r_{0|0}&=&\frac{32(1+k_1b^2)(1+k_2b^2)\big[k_1k_2(3k_2+k_1)b^4+8k_1k_2b^2+3k_1+k_2\big]}{\big[4+(k_1+3k_2)b^2\big]^3}s_0^2\\
 &&+\frac{b^2(3k_1+k_2+4k_1k_2b^2)}{4+(k_1+3k_2)b^2}s_{0|0},\\
 q_{00}&=&\frac{4k_1k_2b^2+3k_1+k_2}{4+(k_1+3k_2)b^2}(\beta
 t_0+s_0^2), \ \ \ q_{0}=\frac{(4k_1k_2b^2+3k_1+k_2)b^2}{4+(k_1+3k_2)b^2}
 t_0.
\eeqn

Finally, assuming (\ref{r00}) holds, we will use the formulas in
Lemma \ref{lem21} to compute the following quantities
$$b^ks_{0|k},\ \ s^m_{\ 0|m}, \ \
s_{0|0}.$$
 Since $n=2$, we always have
  \be\label{y28}
 Ric_{\alpha}=\lambda \alpha^2,
  \ee
where $\lambda=\lambda(x)$ is a scalar function. By (\ref{r00}),
we can get $r^m$ and then we have
  \beq\label{y29}
 r^m_{\ |m}&=&\frac{32(1+k_1b^2)(1+k_2b^2)\big[k_1k_2(3k_2+k_1)b^4+8k_1k_2b^2+3k_1+k_2\big]}{[4+(k_1+3k_2)b^2]^3}s_ms^m\nonumber\\
  &&+\frac{(3k_1+k_2+4k_1k_2b^2)b^2}{4+(k_1+3k_2)b^2}s^m_{\ |m}.
  \eeq
Plugging (\ref{r00}), (\ref{y28}) and (\ref{y29}) into the fourth
formula in Lemma \ref{lem21} and using (\ref{y25}), we obtain
 \beq\label{y30}
s^m_{\ |m}&=&
\frac{8(1+k_1b^2)(1+k_2b^2)[(3k_2^2-k_1^2+6k_1k_2)b^4+4(k_1+3k_2)b^2+8]}{[4+(k_1+3k_2)b^2]^2[2(1-k_1k_2b^4)
+(k_2-k_1)b^2]}\theta\nonumber\\
&&+\frac{b^2[4+(k_1+3k_2)b^2]}{2[2(1-k_1k_2b^4)+(k_2-k_1)b^2]}\lambda.
 \eeq
Plugging (\ref{r00}) into the third formula in Lemma \ref{lem21}
and using (\ref{y25}), we get
 \be\label{y31}
b^ks_{0|k}=\frac{2[2(1-k_1k_2b^4)+(k_2-k_1)b^2]}{4+(k_1+3k_2)b^2}\theta\beta
 \ee
By (\ref{r00}), we can first get $r^m_{\ 0|m}$ expressed in terms
of $t_0,\ s^m_{\ |m},\ s_ms^m$ and $b^ms_{0|m}$. Then plugging
$r^m_{\ m|0}$ ($=0$), $r^m_{\ 0|m}$, (\ref{y30}) and (\ref{y31})
into the second formula in Lemma \ref{lem21} and using
(\ref{y25}), we obtain
 \be\label{y32}
 s^m_{\
 0|m}=\frac{\beta}{2(1-k_1k_2b^4)+(k_2-k_1)b^2]}\Big\{\frac{2A}{[4+(k_1+3k_2)b^2]^2}\theta+\frac{1}{2}[4+(3k_2+k_1)b^2]\lambda\Big\},
 \ee
where
 \beqn
 A:&=&2k_1k_2(k_1^2-18k_1k_2-15k_2^2)b^6+(3k_1^3-135k_1k_2^2-57k_1^2k_2-3k_2^3)b^4\\
 &&-(156k_1k_2+18k_1^2+18k_2^2)b^2-16(3k_1+k_2).
 \eeqn

\begin{lem}
In dimension $n=2$, for a pair $\alpha=\sqrt{a_{ij}y^iy^j}$ and
 $\beta=b_iy^i$, there holds
  \be\label{y33}
 s_{0|0}=b^m(r_{m0|0}-r_{00|m})-\lambda(b^2\alpha^2-\beta^2)+q_{00}-t_{00},
  \ee
  where $\lambda$ is the sectional curvature of $\alpha$.
 \end{lem}

{\it Proof :} By definition, we easily get
 $$s_{0|0}=b^ms_{m0|0}+q_{00}-t_{00}.$$
It follows from Lemma \ref{lem21} that
 $$s_{i0|0}=r_{i0|0}-r_{00|i}-b^m\bar{R}_{0mi0}.$$
Since $n=2$, we have
 $$\bar{R}_{jmik}=\lambda(a_{jk}a_{mi}-a_{ij}a_{mk}).$$
Therefore, we easily obtain (\ref{y33}).  \qed

\

Now by (\ref{r00}), we can first get $r_{m0|0}$ and $r_{00|m}$.
Then plugging them into (\ref{y33}) and using (\ref{y25}) and
(\ref{y31}), we obtain
 \beq\label{y34}
 s_{0|0}&=&\frac{-[4+(k_1+3k_2)b^2]\lambda}{2[2(1-k_1k_2b^4)
+(k_2-k_1)b^2]}(b^2\alpha^2-\beta^2)+2\theta\Big\{\frac{[2(1-k_1k_2b^4)
+(k_2-k_1)b^2]}{4+(k_1+3k_2)b^2}\alpha^2\nonumber\\
&&-\frac{A_1}{[4+(k_1+3k_2)b^2]^2[2(1-k_1k_2b^4)
+(k_2-k_1)b^2]}(b^2\alpha^2-\beta^2)\Big\},
 \eeq
where
 \beqn
 A_1:&=&4k_1k_2(k_1+3k_2)b^6[2k_1k_2b^2+3(k_1-k_2)]+(6k_1^3-6k_2^3-18k_1^2k_2-174k_1k_2^2)b^4\\
 &&-(12k_1^2+28k_2^2+216k_1k_2)b^2-24(3k_1+k_2).
 \eeqn

Now we have obtained the expressions of all those quantities in
(\ref{y024}), as shown in (\ref{y24}), (\ref{y25}), (\ref{y28}),
(\ref{y31}), (\ref{y32}), (\ref{y34}), etc. Plug all those
quantities and (\ref{r00}) into (\ref{K}) and then the flag
curvature {\bf K} of $F$ in Theorem \ref{th1} is given by
 \be\label{y36}
 {\bf
 K}=\frac{2(1+k_2b^2)}{2(1-k_1k_2b^4)+(k_2-k_1)b^2}\Big\{\lambda-\frac{8A_2}{b^2[4+(k_1+3k_2)b^2]^2}s_ms^m\Big\},
 \ee
where
 $$
 A_2:=2k_1^2k_2^2b^6+5k_1k_2(k_1+k_2)b^4+k_1(k_1+13k_2)b^2+2(2k_1+k_2).
 $$

\

In the final step, we will show (\ref{y36}) can be written as
(\ref{y6}). Put $\alpha$ and $\beta$ as
 \be\label{y37}
 \alpha=e^{\sigma}\sqrt{(y^1)^2+(y^2)^2},\  \ \
 \beta=e^{\sigma}(\xi y^1+\eta y^2),
 \ee
 where $\sigma=\sigma(x),\ \xi=\xi(x),\ \eta=\eta(x)$ are scalar
 functions. In the following, we define
 $$
\sigma_{i}:=\sigma_{x^i},\ \ \ \sigma_{ij}:=\sigma_{x^ix^j}, \ \
etc.
 $$
The sectional curvature $\lambda$ of $\alpha$ and the
 norm $b=||\beta||_{\alpha}$ are given by
  \be\label{y38}
 \lambda=-e^{-2\sigma}(\sigma_{11}+\sigma_{22}),\ \ \
 b^2=\xi^2+\eta^2.
  \ee
   Further, we have
   \be\label{y39}
  s_ms^m=\frac{(\xi^2+\eta^2)(\xi_2+\xi\sigma_2-\eta\sigma_1-\eta_1)^2}{4e^{2\sigma}}.
   \ee
 Now comparing (\ref{y37}) with (\ref{y3}) and (\ref{y4}), we get
  \be\label{y40}
 \xi=\frac{u\sqrt{B}}{\sqrt{u^2+v^2}},\ \ \ \eta=\frac{v\sqrt{B}}{\sqrt{u^2+v^2}},
 \ \ \
 \sigma=ln\frac{\sqrt{B}}{\sqrt{u^2+v^2}(1+k_1B)^{\frac{3}{4}}(1+k_2B)^{\frac{1}{4}}}.
  \ee
Plugging  (\ref{y40}) into (\ref{y39}) and using $u_1=v_2,\
u_2=-v_1$, we have
 \be\label{y41}
 s_ms^m=\frac{[4+(k_1+3k_2)B]^2(uB_2-vB_1)^2}{64B\sqrt{1+k_1B}(1+k_2B)^{\frac{3}{2}}}.
 \ee
By (\ref{y5}) we have
 \be\label{y42}
 u_1=v_2,\ \ \ u_2=-v_1,\ \ \ u_{11}+u_{22}=0,\ \ \ v_{11}+v_{22}=0,\ \ \ B_2=-\frac{u}{v}B_1.
 \ee
Plugging $\sigma$ in (\ref{y22}) or (\ref{y40}) into $\lambda$ in
(\ref{y38}) and using (\ref{y42}), we get
 \beq\label{y43}
\lambda&=&-\frac{u^2+v^2}{4B^2\sqrt{1+k_2B}}\Big\{\big[2(1-k_1k_2b^4)
+(k_2-k_1)b^2\big]\sqrt{1+k_1B}(B_{11}+B_{22})\nonumber\\
&&+\frac{(u^2+v^2)B_1^2}{Bv^2(1+k_2B)\sqrt{1+k_1B}}\ T\Big\},
 \eeq
where
$$
 T:=2k_1k_2B^3(k_1k_2B+k_1-k_2)+(k_1^2-k_2^2-8k_1k_2)B^2-4(k_1+k_2)B-2
$$
Now plugging (\ref{y41}) and (\ref{y43}) into (\ref{y36}) and
using (\ref{y42}), we obtain (\ref{y6}).      \qed

\noindent Guojun Yang \\
Department of Mathematics \\
Sichuan University \\
Chengdu 610064, P. R. China \\
 ygjsl2000@yahoo.com.cn


\begin{thebibliography}{999}


\bibitem{BR2} D. Bao and C. Robles, {\it Ricci and flag curvatures
in Finsler geometry}, In "A sampler of Finsler geometry" MSRI
series, Cambridge University Press, 2005.

\bibitem{BRS} D. Bao, C. Robles and Z. Shen, {\it Zermelo
navigation on Riemann manifolds}, J. Diff. Geom. {\bf 66} (2004),
391-449.


\bibitem{CSZ} B. Chen, Z. Shen and L. Zhao, {\it On a calss of Ricci
flat Finsler metrics in Finsler geometry}, preprint.

\bibitem{CS1} X. Cheng and Z. Shen, {\it A class of Finsler metrics with isotropic S-curvature},
Israel Journal of Mathematics, {\bf 169} (2009), 317-340.

\bibitem{CMS} X. Cheng, X. Mo and Z. Shen,
{\it On the flag curvature of Finsler metrics of scalar
curvature}, J. London. Math. Soc., {\bf 68}(2)(2003), 762-780.



\bibitem{CST} X. Cheng, Z. Shen and Y. Tian, {\it A class of
Einstein $(\alpha,\beta)$-metrics}, Isreal J. of Math., {\bf 192}
(1) (2012), 221-249.

\bibitem{shen2}  Z. Shen, {\it  On a class of Landsberg metrics in Finsler geometry},
 Canadian J. of Math., {\bf 61}(6) (2009), 1357-1374.

\bibitem{shen1} Z.  Shen, {\it Volume compsrison and its applications
in Riemann-Finsler geometry}, Adv. in Math. {\bf 128}(1997),
306-328.

\bibitem{shen5} Z. Shen, {\it Finsler metrics with ${\bf K}=0$ and ${\bf S}=0$},
 Canad. J. Math., {\bf 55}(1)(2003), 112-132.


\bibitem{Y} G. Yang, {\it On Randers metrics of isotropic
S-curvature}, Acta Math. Sin., {\bf 52}(6)(2009), 1147-1156 (in
Chinese).

\bibitem{Y1} G. Yang, {\it On Randers metrics of isotropic
S-curvature II}, Publ. Math. Debreceen, {\bf 78}(1) (2011), 71-87.


\bibitem{Y2} G. Yang, {\it A note on a class of   Finsler metrics of isotropic
S-curvature}, preprint.

\bibitem{Y3} G. Yang, {\it On a class of Finsler metrics of
Einstein-reversibility}, preprint.


\bibitem{Y4} G. Yang and X. Cheng, {\it Conformal invariances of two-dimensional Finsler spaces with isotropic main
scalar}, Publ. Math. Debrecen, {\bf 81} (3-4) (2012), 327-340.

\bibitem{Y5} G. Yang, {\it  On a class of two-dimensional Douglas and projectively flat Finsler
metrics}, The Sci. World J., 2013, 291491 [11 pages]
 DOI: 10.1155/2013/291491.

\bibitem{Y6} G. Yang, {\it On a class of two-dimensional singular Douglas and projectively flat Finsler
metrics}, preprint.



















\end{thebibliography}
\end{document}